

\baselineskip=14pt
\parskip=10pt

\font\eighttt=cmtt8
\magnification=\magstephalf

\def\1{{\overline{1}}}
\def\2{{\overline{2}}}
\parindent=0pt
\overfullrule=0in

\def\frac#1#2{{#1 \over #2}}
\centerline
{\bf Automated Proofs of Many Conjectured Recurrences in the OEIS made by R.J. Mathar  }
\bigskip
\centerline
{\it Shalosh B. EKHAD, Mingjia YANG and Doron ZEILBERGER}

The On-Line Encyclopedia Of Integer Sequence (OEIS) ([Sl]), that wonderful resource 
that most combinatorialists, and many other mathematicians and scientists, use
at least once a day, is a treasure trove of mathematical information, and one of its
charms is that it contains many intriguing conjectures. But one should be on one's
guard, because some of the conjectures are either already theorems, or can be routinely proved.

As pointed out in [Z1], all identities involving sequences that belong to the {\bf $C$-finite ansatz}
can be proved by checking sufficiently many initial special cases, and the number of such checks is usually
very small, and the number of initial values needed to make the ``empirical" proof  rigorous can be
easily found.

Recall that a sequence is $C$-finite if it satisfies a homogeneous linear recurrence equation with {\bf constant} coefficients,
the most famous members being sequence $A000079$ \hfill\break
({\tt https://oeis.org/A000079}), namely
$a(n)=2^n$, satisfying the recurrence $a(n)-2a(n-1)=0$ and
sequence $A000045$ ({\tt https://oeis.org/A000045}), satisfying the recurrence $F(n)-F(n-1)-F(n-2)=0$.

Things are not so simple with the $P$-recursive ansatz (see [Z2]), the class of sequences
satisfying a homogeneous linear recurrence equation with {\bf polynomial} coefficients, the most famous
one being sequences A000142 ({\tt https://oeis.org/A000142}), $n!$, satisfying the recurrence
$a(n)\, - \, n\, a(n-1)=0$. Many sequences can be shown, by purely theoretical  {\it a priori} arguments,   to
be $P$-recursive, so it seems that in order to prove that two differently defined $P$-recursive
sequences $a(n)$ and $b(n)$ are actually equal, it should suffice to check sufficiently many
initial values. Alas, while it is true that the sequence $c(n):=a(n)-b(n)$ is also $P$-recursive,
if you don't know the explicit recurrence that it satisfies, it is conceivable (albeit extremely unlikely) that
even if it satisfies a low-order, say, second-order, recurrence, it is something like

$$
( \, n-10^{1000000} \, )\, a(n) \, + \, (n+1) a(n-1) \, + \, (n-3) \, (n-2)=0 \quad ,
$$
in other words it has a `singularity' at a very large positive integer, and the `finitely many' values one has
to check is $10^{1000000}+2$ rather than $2$.

In such cases, to be completely rigorous, one has to actually find the  recurrences, and check that
they do not have such singularities, or if they do, find the largest positive integer that is a singularity.

One case where it is relatively painless (for the computer) to  explicitly find a recurrence for a sequence,
is the case of what is called the {\bf Sch\"utzenberger ansatz} in [Z2], i.e. sequences $\{a(n)\}_0^{\infty}$ whose
{\bf ordinary generating functions}
$$
f(x) = \sum_{n=0}^{\infty} a(n)\,x^n \quad , 
$$
satisfy equations of the form $P(f(x),x)=0$, where $P$ is a polynomial of two variables. More
explicitly, there exists a positive integer $A$, and polynomials $p_i(x)$, $i=0,1,\dots, A$, such that

$$
\sum_{i=0}^A p_i(x) f(x)^i \, = \, 0 \quad .
$$
These are discussed in the Flajolet-Sedgewick {\bf bible}, [FS].

By a well-known algorithm (that by today's standards is fairly straightforward), that goes back to the 19th century, and
is implemented, {\it inter alia}, in Bruno Salvy and Paul Zimmermann's Maple package {\tt gfun}, as
function {\tt algtodiffeq}, every algebraic formal power series satisfies a {\bf linear differential equation}
with {\bf polynomial coefficients}
$$
\sum_{i=0}^L q_i(x) \frac{d^i}{dx^i} f(x) \, = \, 0 \quad ,
$$
that immediately translates (by writing $f(x)=\sum_{n=0}^{\infty} a(n)x^n$, substituting into the
above differential equation, and setting the coefficient of $x^n$ to $0$), to a {\bf linear recurrence equation
with polynomial coefficients}, satisfied by the sequence $a(n)$
$$
\sum_{i=0}^{M} c_i(n) a(n-i) \, = \, 0   \quad .
$$
This is implemented in  the function {\tt diffeqtorec} in {\tt gfun}.

{\bf Mathar's conjectures}

When we searched the OEIS, on July 7, 2017, for

{\tt Conjecture AND recurrence AND  Mathar },

 we got $450$ hits.  Many of them concern the sequences of coefficients of generating functions given in terms of radicals.
While Abel, Galois, and Ruffini famously proved that not every algebraic function can be expressed in terms of radicals,
the converse is trivially true.

Let's take for example, sequence $A004148$ ({\tt https://oeis.org/A004148}),
where a  generating function (conjectured by Michael Somos) is given
$$
f(x)= {\frac {1-x+{x}^{2}-\sqrt {1-2\,x-{x}^{2}-2\,{x}^{3}+{x}^{4}}}{2 \, {x}^{2}}} \quad.
$$

Obviously this is an algebraic formal power series, and the algebraic equation satisfied by $f(x)$ may be routinely 
obtained by clearing radicals. Alas, one often gets a higher order recurrence than the one conjectured.
To prove that the simpler conjectured recurrence also satisfies this sequence
one proceeds as follows.

{\bf Proving that a Lower-Order Recurrence is Equivalent to a Higher-Order Recurrence}

Let $N$ be the forward shift operator:
$$
Na(n):=a(n+1) \quad ,
$$
then the fact that the sequence $\{a(n)\}_0^{\infty}$ satisfies a recurrence
$$
\sum_{i=0}^{M} c_i(n) a(n-i) \, = \, 0   \quad ,
$$
is equivalent to the fact that
$$
\left ( \sum_{i=0}^{M} c_i(n)N^{-i} \right ) a(n) \, \equiv \, 0 \quad .
$$
Calling the operator on the left $C(n,N)$, we have
$$
C(n,N) a(n) \, \equiv \, 0 \quad .
$$

The class of linear recurrence operators with polynomial coefficients is a {\bf non-commutative} algebra, where
everything commutes {\bf except} $n$ and $N$, where the commutation relation between $n$ and $N$ is $N n \,- \, n N \, = \, N$.
See [Z3] for a primer.

We say that the sequence $a(n)$ is {\it annihilated} by the operator $C(n,N)$. Of course, any
{\bf left-multiple} of an annihilator is yet-another annihilator. For example, since
$$
(N-1)(N^2-N-1)=N^3- 2N^2+1 \quad,
$$
the Fibonacci numbers {\it can} be also defined by the recurrence, and initial conditions
$$
F_{n+3}=2F_{n+2}-F_n \quad ; \quad F_{0}=0 \, , \, F_1=1  \, , \, F_2=1 \quad ,
$$
but it would not make William of Ockham happy.

The Euclidean division algorithm for the commutative case is easily  generalized to the case of
the non-commutative algebra of linear recurrence operators with rational-function coefficients.
If  $A(n,N)$ is such a (monic) operator of order $d_1$ and $B(n,N)$ is such a (monic) operator of order $d_2$, with $d_1>d_2$ then there
exist  operators $Q(n,N)$ (the {\it quotient}), of order $d_1 -d_2$, and $R(n,N)$ (of order $<d_2$), 
the {\it remainder}, such that
$$
A(n,N)= Q(n,N) B(n,N)+ R(n,N) \quad .
$$

So suppose that we have a rigorous proof that the operator $A(n,N)$ annihilates our sequence, but by pure guessing,
we found that the sequence is also annihilated by an operator of lower order, $B(n,N)$, in other words, there
exists a lower-order recurrence. If, using this algorithm (that we implemented), it turns out that $R(n,N)=0$, it would follow that
$$
A(n,N)= Q(n,N) B(n,N) \quad .
$$

Since we know, rigorously, that $A(n,N)a(n) \equiv 0$, it follows that
$$
(Q(n,N) B(n,N))a(n) \equiv 0 \quad .
$$
Hence
$$
Q(n,N) ( \, B(n,N) a(n) \, ) \equiv 0 \quad .
$$
Defining $b(n):=B(n,N)a(n)$, we get $Q(n,N) b(n)=0$, and since we already know that $b(i)=0$ for the first few values of $i$,
we have a rigorous proof that $b(n)=0$ for all $n$, and hence that $a(n)$ is annihilated by the lower-order linear
recurrence operator $B(n,N)$.

{\bf The Maple package SCHUTZENBERGER.txt}

While many of the needed functions can be found in the Maple package {\tt gfun} [SZ] mentioned above,
we found it more convenient to write our own Maple code,   {\tt SCHUTZENBERGER.txt}, available from

{\tt http://www.math.rutgers.edu/\~{}zeilberg/tokhniot/SCHUTZENBERGER.txt} \quad .

This is an extension of a Maple package (of the same name) written (many years ago)
by Doron Zeilberger. It contains new procedures to automatically
and effortlessly prove any Mathar-type conjecture. 

{\bf Note:}
In fact, if you are lucky enough to have access to  {\tt Maple Version 12}, we recommend
that you use instead:

{\tt http://www.math.rutgers.edu/~zeilberg/tokhniot/SCHUTZENBERGER12.txt } \quad  ,

that, at least for our purposes, is much better.

The new procedures are listed by typing {\tt ezraOEIS();},
while the old ones can be seen by typing {\tt ezra();}. To get  help on any of these procedures,
with an example, type: \hfill\break {\tt ezra(ProcedureName); } . 

Procedure {\tt MatharConjs()}, is a compilation of $33$ such generating functions for which
R.J. Mathar conjectured recurrences, that readers are welcome to extend.

Procedure  {\tt OEISpaper(L,P,x,n,a)}, inputs such a list {\tt L} (e.g.  {\tt MatharConjs()}) and
it outputs a humanly-readable article, ready for submitting, with statements and
proofs of all the conjectured recurrence listed in $L$. For example, typing:

{\tt OEISpaper(MatharConjs(),P,x,n,a);}

outputs, in about 3 seconds, an article with $33$ theorems, that can be seen in the output file

{\tt http://sites.math.rutgers.edu/\~{}zeilberg/tokhniot/oMathar1maple12.txt} \quad .

[Note that in a few cases the recurrences that we found differ from those of Mathar, but we checked that they are equivalent.]

{\bf An example using the Maple package SCHUTZENBERGER.txt to prove one of Mathar's conjectures}

Have a Maple window, where you have uploaded the package, then type:

{\tt read `SCHUTZENBERGER.txt`}: \quad ,

then open another window, connected to the OEIS. Suppose that you see
that R.J. Mathar conjectured a recurrence, for example OEIS sequence A000957, ({\tt https://oeis.org/A000957}).
Now search for ``G.f'' getting that the (ordinary) generating function for that sequence is 

{\tt (1-sqrt(1-4*x))/(3-sqrt(1-4*x));}

(for our human readers this is  ${\frac {1-\sqrt {1-4\,x}}{3-\sqrt {1-4\,x}}}$).

Then copy-and-paste this generating function into the Maple window, calling it, say, {\tt f};

{\tt f:=(1-sqrt(1-4*x))/(3-sqrt(1-4*x));}

To get the recurrence conjectured by Mathar in {\it computerese}, type:

{\tt radtorec(f,x,n,N);} \quad ,

to get a {\bf statement} of the recurrence, in {\it humanese}, type: 

{\tt  radtorecV(f,P,x,n,N,a):} \quad ,

and to get a {\bf statement and proof} of the recurrence, still in {\it humanese}, type: 

{\tt radtorecVwp(f,P,x,n,a):} \quad  .
\bigskip

{\bf Conclusion} 

While the fact that the sequence of coefficients of any generating function given in terms of radicals
{\bf always} satisfies {\bf some} linear recurrence equation with polynomial coefficients is well known, it is apparently
not as well-known as it should be, or else such recurrences would not be  listed as ``conjectures" in the OEIS.

In the present case, before this article, for each and every such generating function, one had to work pretty
hard to combine all the needed ingredients. After this article (and, more importantly, the new-improved
version of {\tt SCHUTZENBERGER.txt}), it can be done with just one command!

We believe that there exist quite a few other families of ``conjectures" that are routinely provable
and listed in the OEIS (and elsewhere) as ``conjectures". Since
we know that if {\it someone} was ``crazy" enough to actually go through the trouble of proving them
using existing tools, it can be done, {\bf why bother}? Since we know that a proof of the ``conjecture" {\bf exists},
checking it for the first $300$ cases should suffice, and we propose to have a new category called
{\bf ``provable conjecture"}, indicating that we know that there {\bf exists} a proof, but we are
too lazy to actually spell it out, and we are more than happy with  the {\bf empirical verification}.

In other words, the present article is a {\it case study} showing how it is done for this particular family of conjectures.
But enough is enough. We don't see the point of doing it for other families of {\it provable conjectures}. Knowing the fact that it is
provable, plus convincing empirical evidence, suffices! Let's focus on trying to prove conjectures that are not (yet) known to
be provable.

{\bf References}

[FS] P. Flajolet and R. Sedgewick, {\it ``Analytic Combinatorics''}, Cambridge University Press, 2009. \hfill \break
[Freely(!) available on-line from {\eighttt http://algo.inria.fr/flajolet/Publications/book.pdf} ]

[SZ] B. Salvy and P. Zimmermann, {\it Gfun: a Maple package
for the manipulation of generating and holonomic functions
in one variable}, ACM Trans. Math. Software {\bf 20}
(1994).

[Sl] N. J. A. Sloane, {\it The On-Line Encyclopedia of Integer Sequences} (OEIS),  {\tt oeis.org} .

[Z1] Doron Zeilberger, {\it The C-finite Ansatz}, Ramanujan Journal {\bf 31}(2013), 23-32.  \hfill\break
{\tt http://sites.math.rutgers.edu/\~{}zeilberg/mamarim/mamarimhtml/cfinite.html} \quad .

[Z2] Doron Zeilberger, {\it  An Enquiry Concerning Human (and Computer!) [Mathematical] Understanding}, 
in: C.S. Calude, ed., ``Randomness and Complexity, from Leibniz to Chaitin" World Scientific, Singapore, 2007, pp. 383-410. \hfill\break
{\tt http://sites.math.rutgers.edu/\~{}zeilberg/mamarim/mamarimhtml/enquiry.html} \quad .

[Z3] Doron Zeilberger,  {\it Three recitations on holonomic functions and hypergeometric series}, J. Symbolic Comp {\bf 20} (1995), 
699-724 (originally appeared in 24th S\'eminaire Lotharingien, (Spring 1990.) \hfill \break
{\tt http://sites.math.rutgers.edu/\~{}zeilberg/mamarim/mamarimhtml/loth.html} \quad .

\bigskip
\hrule
\bigskip
Shalosh B. Ekhad, c/o D. Zeilberger, Department of Mathematics, Rutgers University (New Brunswick), Hill Center-Busch Campus, 110 Frelinghuysen
Rd., Piscataway, NJ 08854-8019, USA. \hfill\break {\tt ShaloshBEkhad at gmail dot com}   \quad .

\smallskip
Mingjia Yang , Department of Mathematics, Rutgers University (New Brunswick), Hill Center-Busch Campus, 110 Frelinghuysen
Rd., Piscataway, NJ 08854-8019, USA.  \hfill\break Email: {\tt my237 at math dot rutgers dot edu} .

\bigskip
Doron Zeilberger, Department of Mathematics, Rutgers University (New Brunswick), Hill Center-Busch Campus, 110 Frelinghuysen
Rd., Piscataway, NJ 08854-8019, USA. \hfill\break
Email: {\tt DoronZeil at gmail  dot com}   \quad .
\bigskip


\end